\documentclass[a4paper]{jpconf}
\usepackage{graphicx}

\usepackage{amsmath}
\usepackage{latexsym}
\usepackage{amssymb}

\usepackage{amsthm}

\newtheorem{theor}{\theoremname}[section]
\newtheorem{propo}[theor]{\propositionname}
\newtheorem{coro}[theor]{\corollaryname}
\newtheorem{lemm}[theor]{\lemmaname}
\newtheorem{rema}[theor]{\remarkname}

\newenvironment{thm}{\begin{theor}\it}{\end{theor}}

\newenvironment{prop}{\begin{propo}\it}{\end{propo}}

\newenvironment{lem}{\begin{lemm}\it}{\end{lemm}}

\newenvironment{rem}{\begin{rema}\rm}{\end{rema}}

\newtheorem{defin}[theor]{\definitionname} 

\newenvironment{defn}{\begin{defin}\rm}{\end{defin}}

\newtheorem{exa}[theor]{\examplename}

\newenvironment{ex}{\begin{exa}\rm}{\end{exa}}

\newtheorem{exas}[theor]{\examplesname}

\newtheorem{conj}[theor]{\conjecturename}

\newtheorem{pro}[theor]{\problemname}

\newtheorem{prs}[theor]{\problemsname}

\newenvironment{prf}{\begin{proof}}{\end{proof}}
 
\newcommand{\Sesq}{\mathop{{\rm Sesq}}\nolimits}

\renewcommand{\phi}{\varphi}

\font\tengoth=eufm10 at 10pt
\font\sevengoth=eufm7 at 6pt
\newfam\gothfam
\textfont\gothfam=\tengoth
\scriptfont\gothfam=\sevengoth

\newcommand{\mlabel}[1]{\marginpar{#1}\label{#1}}

\newcommand{\fB}{{\mathfrak B}}

\newcommand{\g}{{\mathfrak g}}

\newcommand{\fh}{{\mathfrak h}}

\newcommand{\fq}{{\mathfrak q}}

\renewcommand{\:}{\colon}
\newcommand{\1}{\mathbf{1}}

\newcommand{\cA}{\mathcal{A}}

\newcommand{\cD}{\mathcal{D}}
\newcommand{\cE}{\mathcal{E}}
\newcommand{\cF}{\mathcal{F}}

\newcommand{\cH}{\mathcal{H}}
\newcommand{\cK}{\mathcal{K}}
\newcommand{\cL}{\mathcal{L}}
\newcommand{\cM}{\mathcal{M}}
\newcommand{\cN}{\mathcal{N}}
\newcommand{\cO}{\mathcal{O}}
\newcommand{\cP}{\mathcal{P}}

\newcommand{\cT}{\mathcal{T}}

\newcommand{\cV}{\mathcal{V}}
\newcommand{\cW}{\mathcal{W}}

\newcommand{\dd}{{\tt d}}

\newcommand{\subeq}{\subseteq}

\newcommand{\into}{\hookrightarrow}
\newcommand{\eps}{\varepsilon}

\newcommand{\N}{{\mathbb N}}
\newcommand{\Z}{{\mathbb Z}}
\newcommand{\R}{{\mathbb R}}
\newcommand{\C}{{\mathbb C}}

\newcommand{\T}{{\mathbb T}}

\renewcommand{\hat}{\widehat}

\renewcommand{\tilde}{\widetilde}

\renewcommand{\L}{\mathop{\bf L{}}\nolimits}

\newcommand{\SO}{\mathop{{\rm SO}}\nolimits}

\newcommand{\OO}{\mathop{\rm O{}}\nolimits}

\newcommand{\U}{\mathop{\rm U{}}\nolimits}

\newcommand{\Fix}{\mathop{{\rm Fix}}\nolimits}

\renewcommand{\Re}{\mathop{{\rm Re}}\nolimits}
\renewcommand{\Im}{\mathop{{\rm Im}}\nolimits}

\renewcommand{\tr}{\mathop{{\rm tr}}\nolimits}

\newcommand{\Herm}{\mathop{{\rm Herm}}\nolimits}

\newcommand{\Aut}{\mathop{{\rm Aut}}\nolimits}

\newcommand{\id}{\mathop{{\rm id}}\nolimits}

\renewcommand{\dim}{\mathop{{\rm dim}}\nolimits}

\newcommand{\im}{\mathop{{\rm im}}\nolimits}

\newcommand{\Rarrow}{\Rightarrow}
 
\newcommand{\oline}{\overline}

\newcommand{\la}{\langle}
\newcommand{\ra}{\rangle}

\newcommand{\res}{\vert}

\newcommand{\ssssarr}{\hbox to 15pt{\rightarrowfill}}
\newcommand{\sssarr}{\hbox to 20pt{\rightarrowfill}}
\newcommand{\ssarr}{\hbox to 30pt{\rightarrowfill}}
\newcommand{\sarr}{\hbox to 40pt{\rightarrowfill}}
\newcommand{\arr}{\hbox to 60pt{\rightarrowfill}}
\newcommand{\larr}{\hbox to 60pt{\leftarrowfill}}
\newcommand{\Arr}{\hbox to 80pt{\rightarrowfill}}

\def\theoremname{Theorem}
\def\propositionname{Proposition}
\def\corollaryname{Corollary}
\def\lemmaname{Lemma}
\def\remarkname{Remark}
\def\conjecturename{Conjecture} 

\def\definitionname{Definition}

\def\examplename{Example}
\def\examplesname{Examples}
\def\problemname{Problem}
\def\problemsname{Problems}

\renewcommand{\mlabel}{\label}

\begin{document}
\title{Reflection positivity for the circle group} 

\author{K-H Neeb$^1$ and G Olafsson$^2$} 
\address{$^1$ Department  Mathematik, FAU Erlangen-N\"urnberg, Cauerstrasse 11, 
91058-Erlangen, Germany} 
\address{$^2$ Department of mathematics, Louisiana State University, 
Baton Rouge, LA 70803, USA}
\ead{neeb@math.fau.de, olafsson@math.lsu.edu}

\begin{abstract} In this note we characterize those unitary one-parameter groups 
$(U^c_t)_{t \in \R}$ which admit euclidean realizations in the sense that they are obtained by the 
analytic continuation process corresponding to reflection positivity from a 
unitary representation $U$ of the circle group. 
These are precisely the ones for which 
there exists an anti-unitary involution $J$ commuting with $U^c$. This provides an 
interesting link with the modular data arising in Tomita--Takesaki theory. 
Introducing the  concept of a positive definite function with values in the space of sesquilinear 
forms, we further establish a link between KMS states and reflection positivity on the circle. 
\end{abstract}

\section{Introduction} 

In this note we continue our investigations of the mathematical 
foundations of \textit{reflection positivity}, a basic concept in constructive quantum 
field theory (\cite{GJ81, JOl98, JOl00, JR08,JR07a, JR07b}). 
Originally, reflection positivity, also called Osterwalder--Schrader positivity, 
arises as a requirement on the euclidean side to establish a 
duality between euclidean and relativistic quantum field theories (\cite{OS73}). 
It is closely related to ``Wick rotations'' or 
``analytic  continuation'' in the time variable 
from the real  to the imaginary axis. 

The underlying fundamental concept is that of a 
{\it reflection positive Hilbert space}, introduced in \cite{NO14a}. 
This is a triple $(\cE,\cE_+,\theta)$, 
where $\cE$ is a Hilbert space, $\theta : \cE \to \cE$ is a unitary involution
and $\cE_+$ is a closed subspace of $\cE$ which is $\theta$-positive in the sense that 
the hermitian form $\langle \theta u,v\rangle$ is
positive semidefinite on $\cE_+$.
\begin{footnote}{We write $\la v, w \ra$ for the scalar product in a 
complex Hilbert space and we assume that it is linear in the {\it first} 
argument.}  
\end{footnote}
Let 
$\cN := \{ v \in \cE_+ \: \la \theta v, v \ra = 0\}$, 
write $\hat\cE$ for the Hilbert space completion of the quotient $\cE_+/\cN$ 
and $q \: \cE_+ \to \hat\cE, v \mapsto \hat v$ for the canonical map. 
If $T \: \cD(T) \subeq \cE_+ \to \cE_+$ 
is a linear or antilinear operator with $T(\cN\cap \cD(T)) \subeq \cN$, 
then there exists a well-defined operator $\hat T \: \cD(\hat T) \subeq \hat\cE \to \hat\cE$ defined by 
$\hat T(\hat v) = \hat{Tv}$ for $v \in \cD(T)$. The\ $\ \hat{}\ $\ is the 
{\it OS-quantization functor.}

To see how this relates to group representations, let us call a triple 
$(G,H,\tau)$ a {\it symmetric Lie group} if $G$ is a Lie group, 
$\tau$ is an involutive automorphism of $G$ and $H$ is an open subgroup 
of the group $G^\tau$ of $\tau$-fixed points. Then the Lie algebra 
$\g$ of $G$ decomposes into $\tau$-eigenspaces $\g = \fh \oplus \fq$ 
and we obtain the {\it Cartan dual Lie algebra} $\g^c=\fh\oplus i\fq$. 
If $(G,H,\tau)$ is a symmetric Lie group and
$(\cE,\cE_+,\theta)$ a reflection positive Hilbert space, then we say that 
a unitary representation $\pi \: G \to \cE$ is {\it reflection positive 
with respect to $(G,H,\tau)$} if the following three conditions hold:
\begin{description}
\item[\rm(RP1)] $\pi(\tau(g)) = \theta \pi(g)\theta$ for every $g \in G$. 
\item[\rm(RP2)] $\pi (h)\cE_+= \cE_+$ for every $h \in H$. 
\item[\rm(RP3)] There exists a subspace $\cD\subeq \cE_+\cap \cE^\infty$, dense in $\cE_+$, 
such that  $\dd \pi (X)\cD\subset \cD$ for all $X\in \fq$. 
\end{description}
A typical source of reflection positive representations are the representations 
$(\pi_\phi, \cH_\phi)$ obtained via GNS construction (cf.~Theorem~\ref{prop:gns}) 
from $\tau$-invariant positive definite functions \break 
$\phi \: G \to B(\cV)$, respectively the kernel $K(x,y) = \phi(xy^{-1})$, 
where $\cV$ is a Hilbert space (cf.\ \ref{sec:a}). If $G_+ \subeq G$ is an open subset 
with $G_+ H = G_+$, then $\phi$ is called {\it reflection positive for $(G,G_+,\tau)$} 
if the kernel 
$Q(x,y) = \phi(x\tau(y)^{-1})$ is positive definite on $G_+$. For 
$\cE = \cH_\phi$, the subspace 
$\cE_+$ is generated by the functions $K(\cdot, y)v$, $y \in G_+, v \in \cV$, 
and $\hat\cE$ identifies naturally with the Hilbert space $\cH_Q \subeq \cV^{G_+}$ 
(cf.\ \cite[Prop.1.11]{NO14a}, \cite[Ex.~5.17a]{MNO14}). If the kernels 
$\la Q(x,y)v,v\ra$ are smooth for $v$ in a dense subspace of $\cV$, then (RP1-3) 
are readily verified. 

If $\pi$ is a reflection positive representation, 
then $\pi^c_H(h):=\hat{\pi (h)}$ defines a unitary representation 
of $H$ on $\hat\cE$. However, we would like to have a 
unitary representation $\pi^c$ of the simply connected Lie group $G^c$ with 
Lie algebra $\g^c$ on $\hat\cE$ extending $\pi^c_H$ in such a way that the 
derived representation is compatible with the operators 
$i\hat{\dd\pi(X)}$, $X \in \fq$, that we obtain from (RP3) on a dense 
subspace of~$\hat\cE$. If such a representation 
exists, then we call 
$(\pi, \cE)$ a {\it euclidean realization} of the representation $(\pi^c,\hat\cE)$ of $G^c$. 
Sufficient conditions for the existence of $\pi^c$ have been developed in 
\cite{MNO14}. The prototypical pair $(G,G^c)$ consists of the euclidean 
motion group $\R^d \rtimes \OO_d(\R)$ and the simply connected covering of the 
Poincar\'e group $\cP_+^{\uparrow} = \R^d \rtimes \SO_{1,d}(\R)_0$. 

In \cite{NO14b} 
we studied reflection positive one-parameter groups and hermitian contractive 
semigroups as one key to reflection positivity for more general 
symmetric Lie groups and their representations. Here a crucial point is that, for 
every unitary one-parameter group $U^c_t = e^{itH}$ with $H \geq 0$ on the Hilbert space 
$\cV$, we obtain by 
$\phi(t) := e^{-|t|H}$ a $B(\cV)$-valued function on $\R$ which is reflection positive 
for $(\R,\R_+, -\id_\R)$ and which leads to a euclidean realization of $U^c$. 
From this we derive that all representations of the $ax+b$-group, resp., the Heisenberg group 
which satisfy the positive spectrum condition for the translation group, resp., 
the center, possess natural euclidean realizations. 

The present note grew out of the attempt to extract the representation theoretic 
aspects from the discussion of reflection positivity for the circle group 
$\T$ in \cite{KL81}. 
Here we continue our project by exploring reflection positive functions  
$\phi \: \T \to B(\cV)$ for $(\T,\T^+, \tau)$, where $\tau(z) = z^{-1}$ and $\T^+$ is a half circle. 
This leads us naturally to anti-unitary involutions, an 
aspect that did not show up for triples $(G,G_+, \tau)$, where $G_+$ is a semigroup. 
We start in Section~\ref{sec:2} with a discussion of the 
integral representation of reflection positive operator-valued functions for $(\T,\T^+, \tau)$ 
due to Klein and Landau (\cite{KL81}) 
and in Section~\ref{sec:3} we characterize those unitary one-parameter groups 
$(U^c,\cH)$ which admit euclidean realizations in this context as those for which 
there exists an anti-unitary involution $J$ commuting with $U^c$ (Theorem~\ref{thm:realize}). 
Any such pair $(J,U^c)$ with $U^c_t = e^{itH}$ can be encoded in the pair $(J,\Delta)$, 
where $\Delta = e^{-\beta H}$ is positive selfadjoint with $J\Delta J = \Delta^{-1}$, a relation 
well-known from Tomita-Takesaki theory. Such 
pairs are closely linked to real standard subspaces, a connection discussed in Section~\ref{sec:4}. 
In Section~\ref{sec:5} we finally explain how all this connects to KMS states for 
$C^*$-dynamical systems $(\cA,\R, \alpha)$ and conclude with a short discussion of perspectives. 
To establish the connection with KMS states, we need the concept of a positive definite 
function on a group $G$ with values in the space $\Sesq(\cV)$ of sesquilinear forms on a 
real or complex vector space $\cV$. This concept is briefly developed in \ref{sec:a}, 
where we explain in particular how the GNS construction works in this context. 
\ref{sec:b} provides some tools to obtain integral representations of 
positive definite functions on convex sets which sharpen some known results in this context. 

We hope that this note will prove useful in the further development of the representation 
theoretic side of reflection positivity under the presence of anti-unitary 
involutions which occur in many recent constructions in Quantum Field Theory 
(see in particular \cite[p.~627]{BS04},  \cite{Bo92}, \cite{BGL02}). 



\section{Reflection positivity on the circle group $\T$} 
\mlabel{sec:2}

Let $G := \T_\beta := \R/\beta \Z$, 
where $\beta > 0$, so that $\T_\beta$ is a circle of length~$\beta$. 
We write $[t] := t + \beta \Z$ 
for the image of $t$ in $\T_\beta$ and write 
$\T_\beta^+ := \{ [t] \in \T_\beta \: 0 < t < \beta/2\}$ 
for the corresponding semicircle. 
We further fix the involutive automorphism $\tau_\beta(z) = z^{-1}$ given by 
inversion in $\T_\beta$. 
In the following we identify functions on $\T_\beta$ with $\beta$-periodic functions on $\R$. 

\begin{defn}
Let $\cV$ be a Hilbert space. 
A weak operator continuous function $\phi \: \T_\beta \to B(\cV)$ is called 
{\it reflection positive} 
w.r.t.\ $(\T_\beta, \T_\beta^+,\tau_\beta)$ 
if and only if it is positive definite and the kernel 
$(\phi(t + s))_{0 < t,s < \beta/2}$ 
is positive definite, which is equivalent to the positive definiteness of the kernel 
\begin{equation} 
  \label{eq:q-phi}
 Q_\phi(t,s) := \phi\Big(\frac{t + s}{2}\Big), \quad  0 < t,s < \beta \quad 
\mbox{ on } ]0,\beta[. 
\end{equation}
\end{defn}

\begin{rem} \mlabel{rem:2.2}
The reflection positivity of $\phi$ 
ensures that the corresponding $\beta$-periodic GNS-representation of $\R$ on the reproducing 
kernel Hilbert space $\cE := \cH_\phi \subeq \cV^{\T_\beta} \subeq \cV^\R$ by 
$\pi_\phi(t)f:= f(\cdot +t)$ is reflection positive with respect to 
$(\theta f)(x) = f(-x)$ and the closed subspace $\cE_+$ generated by 
$\phi(-t)v$, $t \in \T_\beta^+$, $v \in \cV$ 
(cf.\ Proposition~\ref{prop:gns}). Here $\tau_\beta(x) = -x$ and $H = \{0\}$. 
Conditions (RP1/2) are obvious, and (RP3) is satisfied 
because the restriction of $\phi$ to $]0,\beta[$ is automatically smooth 
(Theorem~\ref{thm:2.4}). For the basis element $X= 1 \in \L(\T_\beta) \cong \R$, 
the operator $\hat{\dd\pi(X)}$ acts on $\hat\cE \subeq \cV^{\T_\beta^+}$ by $\frac{d}{dt}$. 
\end{rem} 

\begin{ex} Basic examples of reflection positive 
functions are given by  
\[ f_\lambda(t) = e^{-t\lambda} + e^{-(\beta - t)\lambda} 
= 2 e^{-\beta \lambda/2} \cosh((\textstyle{\frac{\beta}{2}}-t)\lambda)\] 
for $0 \leq t \leq \beta$  and $\lambda \geq 0$. 
The corresponding kernel $Q_{f_\lambda}$ is positive definite 
on all of $\R$ because $f_\lambda$ is a Laplace transform of the positive measure 
$\delta_\lambda + e^{-\beta \lambda} \delta_{-\lambda}$. 
A direct calculation shows that the 
Fourier series of the $\beta$-periodic extension  of 
$f_\lambda$ to $\R$ (also denoted $f_\lambda$) is given by 
\begin{equation}
  \label{eq:fourexp}
f_\lambda(t) = \sum_{n \in\Z} c_n e^{2\pi i nt/\beta} \quad \mbox{ with }  \quad 
c_n = \frac{2\beta \lambda(1 - e^{-\beta \lambda})}{(\lambda\beta)^2 + (2\pi n)^2}.
\end{equation}
As $c_n \geq 0$ for every $n \in \Z$, the function $f_\lambda$ on $\T_\beta$ 
is positive definite. This shows that $f_\lambda$ is reflection 
positive.
\end{ex} 

\begin{thm} {\rm(\cite[Thm.~3.3]{KL81})}  \mlabel{thm:2.4}
A $\beta$-periodic weak operator continuous function \break 
$\phi \: \R \to B(\cV)$ 
is reflection positive w.r.t.\ $(\T_\beta,\T_\beta^+,\tau_\beta)$ 
if and only if there exists a $\Herm(\cV)_+$-valued measure 
$\mu_+$ on $[0,\infty[$ such that 
\begin{equation}
  \label{eq:phi-intrep}
\phi(t) = \int_0^\infty e^{-t\lambda} + e^{-(\beta - t)\lambda}\, d\mu_+(\lambda) 
\quad \mbox{ for } \quad 0 \leq t \leq \beta.  
\end{equation}
Then the measure $\mu_+$ is uniquely determined by $\phi$. 
\end{thm}

\begin{prf} If $\phi$ is given by \eqref{eq:phi-intrep}, 
then the kernel $\phi\big(\frac{x+y}{2}\big)$ is clearly positive definite. 
To see that $\phi$ is a positive definite function of $\T_\beta$, 
we use the Fourier 
expansion $f_\lambda(t) = \sum_n c_n(\lambda) e^{2\pi i nt/\beta}$ from 
\eqref{eq:fourexp} to obtain 
\[ \phi(t) 
= \sum_{n \in \Z}e^{2\pi i nt/\beta}  \int_0^\infty c_n(\lambda)\, d\mu_+(\lambda).\] 
Now the positivity of the operators $\int_0^\infty c_n(\lambda)\, d\mu_+(\lambda)$ 
shows that $\phi$ is positive definite, hence reflection positive. 

If, conversely, $\phi \: [0,\beta] \to B(\cV)$ 
is reflection positive, then it can be written as a 
Laplace transform $\phi = \cL(\mu)$ of a $\Herm(\cV)_+$-valued measure  $\mu$ on $\R$  
(\cite[Thm.~18.8]{Gl03}; see also Theorem~\ref{thm:I.7} below). 
Since the reflection in $\beta/2$ is the composition of the reflection 
$r(t) =-t$ in $0$ and the translation by $\beta$, the function $\phi$   
is also symmetric with respect to $\frac{\beta}{2}$.
This symmetry requirement is equivalent to 
$r_*\mu = e_{-\beta} \mu$ for $e_\beta(\lambda) =e^{\beta \lambda}$, 
so that $\phi$ can be written as in \eqref{eq:phi-intrep}. 
\end{prf}

\begin{rem} \mlabel{rem:meas} 
(a) It is often more convenient to work 
with the measure $\mu$ on $\R$ defined by 
\begin{equation} \label{eq:sum-mu}
d\mu(\lambda) = d\mu_+(\lambda) + e^{\beta\lambda} d\mu_+(-\lambda),  
\quad \mbox{ which satisfies } \quad 
r_*\mu = e_{-\beta} \mu \quad \mbox{ for } \quad r(\lambda) = -\lambda.
\end{equation}
We thus obtain the description 
$\phi(t) = \int_\R e^{-t\lambda}\, d\mu(\lambda) = \cL(\mu)(t)$ 
of $\phi$ as the Laplace transform $\cL(\mu)$. 
The existence of these integrals for $0 \leq t \leq \beta$ 
only requires that the measure 
$\mu$ is finite. We then have $\mu(\R) = \phi(0) = \phi(\beta)$ 
by \eqref{eq:sum-mu}. 

(b) In view of \eqref{eq:sum-mu}, the measure $\nu := e_{-\beta/2}\mu$ is symmetric. 
For the Fourier transform $\hat\mu(z) = \int_\R e^{iz\lambda}\, d\mu(\lambda)$ which is defined on 
the closed strip $\oline{\cD_\beta} = \{ z \in \C \: 0 \leq \Im z \leq \beta\}$, we 
thus obtain  the relations 
\begin{equation}
  \label{eq:fourier}
\hat\mu(i \beta-z) = \hat\mu(z), \quad z \in \oline{\cD_\beta},\quad  
\quad \mbox{ and } \quad \hat\nu(t) = \hat\nu(-t)= \hat\mu\Big(t + i\frac{\beta}{2}\Big), \quad t \in \R.
\end{equation}
%
\end{rem}

\begin{rem} \mlabel{rem:2.6}  (The representation $U^c$ on $\hat\cE$) 
With the integral representation from Theorem~\ref{thm:2.4}, we can make 
the unitary one-parameter group $U^c$ of the dual group $(\T_\beta)^c \cong \R$ on $\hat\cE$ more explicit. 
From Remark~\ref{rem:2.2}, we recall that 
$\hat\cE$ can be identified with the reproducing kernel Hilbert space corresponding 
to the kernel $\phi\big(\frac{x+y}{2}\big)$ on the real interval 
$[0,\beta]$ (cf.\ \eqref{eq:q-phi}).

(a) In view of Theorem~\ref{thm:I.7}, $\hat\cE$ can be identified 
with the vector-valued $L^2$-space $L^2(\R,\mu;\cV)$ with the scalar product 
$\la \xi, \eta \ra = \int_\R \la d\mu(\lambda) \xi(\lambda), \eta(\lambda)\ra.$
On this space we have a natural unitary one-parameter group defined by 
\[ (U^c_t f)(\lambda) := e^{it \lambda} f(\lambda) \quad \mbox{ for }\quad 
t, \lambda \in \R.\]
Since $L^2(\R,\mu;\cV)$ contains the constant functions, we obtain a bounded 
operator $j \: \cV \to L^2(\R,\mu;\cV), j(v)(\lambda) := v$. Then $j(\cV)$ is $U^c$-cyclic 
in $L^2(\R,\mu;\cV)$ (\cite[Lemma~III.8]{Ne98}), 
and we have 
\begin{equation}
  \label{eq:a}
 \la U^c_t j(v), j(w) \ra 
= \int_\R e^{it\lambda} \la d\mu(\lambda)v, w \ra = \la \phi(-it) v,w \ra 
\quad \mbox{ for } \quad v,w \in \cV, t \in \R, 
\end{equation}
where we have used the analytic continuation of $\phi$ to the strip 
$\{ z \in \C \: 0 \leq \Re z \leq \beta\}$ that follows from the 
integral representation (Theorem~\ref{thm:2.4}) and Theorem~\ref{thm:I.7}. 
Therefore $(U^c,L^2(\R,\mu;\cV))$ can be identified with the vector-valued 
GNS representation 
corresponding to the positive definite function $\hat\mu(t) = \phi(-it)$ on $\R$ 
obtained by analytic continuation from $\phi$ (cf.\ Proposition~\ref{prop:gns}). 
If $(Hf)(\lambda) = \lambda f(\lambda)$ is the infinitesimal generator of $U^c_t = e^{itH}$, 
then, for every $v \in \cV$, we have 
$v \in \cD(e^{-\beta H/2})$ because 
$\int_\R e^{-\beta \lambda}\la d\mu(\lambda)v,v\ra = \la \phi(\beta)v,v \ra$ is finite 
(Lemma~\ref{lem:a.2}). 
For $0 \leq t \leq \beta$ and $v,w \in \cV$, we now obtain 
$\la \phi(t)v,w\ra = \int_\R e^{-t\lambda}\la d\mu(\lambda)v,w\ra  
= \la e^{-tH}v,w \ra,$ 
resulting in the dilation formula 
\begin{equation}
  \label{eq:dil}
\phi(t) = j^* e^{-tH} j \quad \mbox{ for }\quad 0 \leq t \leq \beta.
\end{equation}

Further, $(Jf)(\lambda) := e^{-\beta\lambda/2} f(-\lambda)$, is a unitary involution 
on $L^2(\R,\mu;\cV)$ with $JHJ = -H$, and 
\[ R := e^{\beta H/2} J = J e^{-\beta H/2}, \quad (Rf)(\lambda) = f(-\lambda) \] 
is an involution with domain $\cD(R) = \cD(e^{-\beta H/2})$ and $\cV \subeq \Fix(R) :=\{ v \in \cD(R) \: Rv =v\}$. 

(b) Alternatively, we can use $\tilde j \: \cV \to L^2(\R,\mu;\cV), \tilde j(v)(\lambda) = e^{-\beta \lambda/4}v$ 
to obtain 
\[ \la U^c_t \tilde  j(v), \tilde j(w) \ra 
= \int_\R e^{it\lambda} e^{-\beta \lambda/2}\la d\mu(\lambda)v, w \ra 
= \Big\la \phi\Big(\frac{\beta}{2}-it\Big) v,w \Big\ra,\] 
which exhibits $U^c$ as the GNS representation of the symmetric positive definite function 
$\psi(t) := \phi\big(\frac{\beta}{2}-it\big).$
\end{rem}

\section{Existence of euclidean realizations} 
\mlabel{sec:3}

The following proposition provides various characterizations of unitary one-parameter groups 
with reflection symmetry. As we shall see below, these are precisely the ones with a 
euclidean realization from $(\T_\beta,\T_\beta^+, \tau_\beta)$. 

\begin{prop} \mlabel{prop:2.9} 
For a unitary one-parameter group 
$(U_t)_{t \in \R}$ on $\cH$ with spectral measure \break $E \: \fB(\R) \to B(\cH)$, 
the following are equivalent:
\begin{itemize}
\item[\rm(i)] There exists an anti-unitary involution $J$ on $\cH$ with 
$J U_t J = U_t$ for $t \in \R$. 
\item[\rm(ii)] For $\cH_\pm := E(\R^\times_\pm)\cH$, the unitary one-parameter groups 
$U^+_t := U_t\res_{\cH_+}$ and $U^-_t := U_{-t}\res_{\cH_-}$ are unitarily equivalent. 
\item[\rm(iii)] The unitary one-parameter group 
$(U,\cH)$ is equivalent to a GNS representation 
$(\pi_\psi, \cH_\psi)$, where $\psi \: \R \to B(\cV)$ 
is a symmetric positive definite function. \
\item[\rm(iv)] There exists a unitary involution $R$ on $\cH$ with 
$R U_t R = U_{-t}$ for $t \in \R$. 
\end{itemize}
\end{prop}

\begin{prf} Every cyclic unitary one-parameter group is isomorphic to 
$L^2(\R,\mu)$ with $(U_t f)(\lambda) = e^{it\lambda} f(\lambda)$. Then  $Kf = \oline f$ 
is an anti-unitary involution. 
Decomposing into cyclic subspaces, we thus obtain an anti-unitary involution 
$K$ on $\cH$ satisfying $KU_t K = U_{-t}$ for $t \in \R$. 
As (i)-(iv) hold for the trivial 
representation on the subspace $\cH_0 = E(\{0\})\cH$ 
of fixed points, we may w.l.o.g.\ assume that $\cH_0 = \{0\}$, so that 
$\cH = \cH_+ \oplus \cH_-$.  

(i) $\Rarrow$ (ii): The operator $W := KJ$ is unitary and satisfies 
$W U_t W^{-1} = U_{-t}$ for $t \in \R$. 
Therefore $W\cH_\pm = \cH_\mp$ and we obtain a unitary intertwining operator 
$(U^+,\cH_+) \to (U^-,\cH_-)$. 

(ii) $\Rarrow$ (iii): We write $\cH = \cK \oplus \cK$ with 
$U_t = W_t \oplus W_{-t}$ and $(W,\cK) \cong (U^+, \cH_+)$. 
Then the representations on both summands are disjoint. 
Put $\cV := \{ (v,v) \: v \in \cK\}$ and let $\cW \subeq \cH$ be the 
closed $U$-invariant subspace generated by $\cV$. 
Since $\cW$ is invariant under all spectral projections, 
it contains the projection of $\cV$ onto both factors, so that 
$\cW = \cH$. The corresponding positive definite function $\psi(t) := P_\cV U_t P_\cV^*$ 
obtained from the orthogonal projection $P_\cV \: \cH \to \cV$ satisfies 
\[ \la \psi(t)(v,v), (w,w) \ra 
=  \la (W_tv,W_{-t} v), (w,w) \ra 
=  \la (W_t + W_{-t})v,w \ra \quad \mbox{ for }\quad v,w \in \cK, \]
so that $\psi$ is symmetric. Hence (iii) follows from Proposition~\ref{prop:gns}. 

(iii) $\Rarrow$ (iv): If $\psi$ is symmetric, then 
$R \: \cH_\psi \to \cH_\psi,  (Rf)(t) := f(-t)$ 
is a unitary involution with the required properties because the kernel 
$Q(t,s) = \psi(t-s)$ satisfies $Q(-t,-s) = Q(t,s)$ (cf.\ \cite[Rem.~II.4.5(c)]{Ne00}). 

(iv) $\Rarrow$ (ii) follows from $R(\cH_\pm) = \cH_\mp$. 

(ii)  $\Rarrow$ (i): As above, we write $\cH = \cK \oplus \cK$ with 
$U_t = W_t \oplus W_{-t}$ and $(W,\cK) \cong (U^+, \cH_+)$ 
and let $K \: \cK \to \cK$ be an anti-unitary involution with 
$K U_t K = U_{-t}$ for $t \in \R$. Then $J(v,w) := (Kw,Kv)$ has the required 
properties. 
\end{prf}


\begin{rem} Let $U_t = e^{itH}$ be a continuous unitary one-parameter 
group on the complex Hilbert space 
$\cH$ with infinitesimal generator $H$ and $J$ an anti-unitary involution with 
$J U_t J = U_{\eps t}$ for  $\eps \in \{\pm 1\}$. Then $JHJ = - \eps H.$ 
If the operator $H$ is non-negative, then so is $JHJ$, which can only happen for 
$\eps = -1$. 
\end{rem}

\begin{defn}\mlabel{def:1.13} We have seen in Remark~\ref{rem:2.2} that every 
weakly operator continuous reflection positive function 
$\phi \: \T_\beta \to B(\cV)$ leads to a reflection positive unitary representation 
$(U,\cE)$ of $\T_\beta$ (Remark~\ref{rem:2.2}). 
If $(U^c, \hat\cE)$ is the corresponding unitary one-parameter group from 
Remark~\ref{rem:2.6}, then we call 
$(U,\cE)$ a {\it euclidean realization of $U^c$}. 
We have already seen in Remark~\ref{rem:2.6} that it 
can be obtained by the GNS construction from the positive definite function 
$\R \to B(\cV), t \mapsto \phi(it)$ or the symmetric function 
$\phi(it + \frac{\beta}{2})$. In this sense it is obtained from 
$U$ by an analytic continuation process. 
\end{defn}

At this point it is a natural question which unitary one-parameter groups 
$(U^c, \cH)$ have a euclidean realization in the sense of Definition~\ref{def:1.13}. 
This can now be stated 
in terms of the conditions discussed in  Proposition~\ref{prop:2.9}: 

\begin{thm} \mlabel{thm:realize} {\rm(Realization Theorem)} 
A unitary one-parameter group $(U_t^c)_{t \in \R}$ on a Hilbert space 
$\cH$ has a euclidean realization  in terms of a reflection positive 
representation of $(\T_\beta, \T_\beta^+,\theta)$ 
if and only if there exists an anti-unitary 
involution $J$ on $\cH$ commuting with $U^c$. 
\end{thm}

\begin{prf} Since the assertion is trivial for the trivial representation, 
we may assume that there are no non-zero fixed vectors, i.e., $\cH^{U^c} = \{0\}$. 

From Proposition~\ref{prop:2.9} we know that the existence of an 
anti-unitary involution commuting with $U^c$ is equivalent to the realizability of 
$U^c$ by the GNS construction from a symmetric positive definite function. 
If $U^c$ has a euclidean realization as in Definition~\ref{def:1.13}, 
then $\psi(t) := \phi\big(\frac{\beta}{2} + it)$ 
is such a function (Remark~\ref{rem:2.6}(b)). 

Suppose, conversely, that there exists an anti-unitary 
involution $J$ on $\cH$ commuting with $U^c$. 
As in the proof of  Proposition~\ref{prop:2.9}, we 
write $\cH = \cV \oplus \cV$ with $U_t^c = e^{itA} \oplus e^{-itA}$ for a 
positive selfadjoint operator $A$ on $\cV$. For $\beta > 0$, we consider the bounded operator 
$j \: \cV \to \cH, j(v) = (v,e^{-\frac{\beta}{2}A}v).$ 
Since the projection 
$P(v_1, v_2) = (v_1, 0)$ onto the first component is $E(\R^\times_+)$ for the 
spectral measure $E$ of $U^c$ (here we use that $E(\{0\})=0$), the closed 
$U^c$-invariant subspace generated by $j(\cV)$ is adapted to the decomposition 
$\cH = \cV \oplus \cV$. Hence its cyclicity follows from the fact that the range of 
$e^{-\frac{\beta}{2}A}$ is dense in $\cV$. As $j(\cV)$ is cyclic, 
\[ \psi \: \R\to B(\cV), \quad \psi(t) := j^* U^c_t j = j^* e^{itH} j \quad \mbox{ for } \quad 
U^c_t = e^{itH}, \] 
is a positive definite function for which the GNS representation 
is equivalent to $(U^c, \cH)$ (Proposition~\ref{prop:gns}). 
Now $j^*(w_1, w_2) =  w_1 + e^{-\frac{\beta}{2}A}w_2$ yields 
\[ \psi(t) 
= e^{itA} + e^{-\frac{\beta}{2}A} e^{-itA} e^{-\frac{\beta}{2}A} 
= e^{itA} + e^{(-it-\beta)A}.\] 
Since $A \geq 0$, 
$\phi(z) := \psi(iz) = e^{-zA} + e^{-(\beta -z)A}$ 
defines a bounded operator for $0 \leq \Re z \leq \beta$. 
For $0 \leq t \leq \beta$, we thus obtain the positive definite function 
\begin{equation}
  \label{eq:phi1}
\phi(t) = e^{-tA} + e^{-(\beta-t)A} = j^* e^{-tH} j, \quad 0 \leq t \leq \beta,  \quad \mbox{ satisfying } \quad \phi(\beta - t) = \phi(t).
\end{equation}
Hence $\phi$ defines a reflection positive function on $\T_\beta$ for which 
$\phi(-it) = \psi(t)$, so that by \eqref{eq:a} the corresponding 
representation of the dual group $(\T_\beta)^c\cong \R$ is equivalent to~$U^c$. 
\end{prf}

\begin{rem} \mlabel{rem:1.16} 
(a) In \cite[Prop.~6.1]{NO14b}, we have shown that a unitary one-parameter group 
$(U_t)_{t \in \R}$ with $U_t = e^{itH}$ has a euclidean realization in terms of 
the triple $(\R,\R_+, -\id_\R)$ if and only if $H \geq 0$, and then 
$\phi(t) := e^{-|t|H}$ is a corresponding reflection positive function. 

(b) For the function $\phi$ from 
\eqref{eq:phi1}, a $\Herm(\cV)_+$-valued measure with $\cL(\mu) = \phi$ 
is obtained as in \eqref{eq:sum-mu} from the spectral measure $E_+$ of $A$ by 
$d\mu(\lambda) = dE_+(\lambda) + e^{\beta \lambda} dE_+(-\lambda)$. 
Here $e^{-tH}= e^{-tA} \oplus e^{tA}$, so that 
$j(v) = (v,e^{-\beta A/2}v) \in \cD(e^{\beta H/2})$. 

(c)  Let $U^c_t = e^{itH}$  be a unitary one-parameter group on $\cH$ and  
$J$ be a unitary involution on $\cH$ with $J H J = -H$. 
Then $R := J e^{-\beta H/2}$ is an unbounded involution with 
$\cD(R) = \cD(e^{-\beta H/2})$. We further assume that 
$j \: \cV \to \cH$ is a bounded operator with cyclic range contained in 
$\Fix(R) = \{ v \in \cD(R) \: Rv =v\}$. Then 
\begin{equation}
  \label{eq:dil2}
\phi(t) := j^* e^{-tH} j \in B(\cV), \quad 0  \leq t \leq \beta 
\end{equation}
defines a function on $[0,\beta]$ 
for which the kernel $\phi\big(\frac{t+s}{2}\big)$ is positive 
definite, and we further have 
\begin{align*}
\phi(t) 
&= j^* e^{-tH} j 
= j^* R^* e^{-tH} R j 
= j^* e^{-\beta H/2} J e^{-tH} J e^{-\beta H/2} j \\
&= j^* e^{-\beta H/2} e^{tH} e^{-\beta H/2} j 
= j^* e^{-(\beta - t)H} j = \phi(\beta -t).
\end{align*}
Therefore $\phi$ defines a reflection positive function on $\T_\beta$. In view of 
\eqref{eq:dil2}, it is the 
Laplace transform of the measure $j^* E j$, where $E$ is the spectral measure of $H$. 
\end{rem}

\section{Standard subspaces} 
\mlabel{sec:4}

To connect reflection positivity on $\T_\beta$ with the modular theory of von Neumann 
algebras, we now take a closer look at {\it standard real subspaces} of a complex Hilbert space $\cH$. 

\begin{defn} A closed real subspace $V \subeq \cH$ is said to be {\it standard} if 
$V \cap i V = \{0\}$ and $\oline{V + i V} = \cH.$ Then we define the corresponding 
{\it Tomita operator} on the dense subspace $V_\C := V + i V$ of $\cH$ 
by $S_V(x + iy) := x-iy$ for $x,y \in V$ (cf.~\cite{BGL02}). 
\end{defn}

\begin{lem} \mlabel{lem:e.1} If $V \subeq \cH$ is a standard real subspace and 
$S := S_V$, then 
\begin{itemize}
\item[\rm(i)] $S$ is a closed densely defined antilinear involution. 
\item[\rm(ii)] If $\Delta := S^*S$ and $S = J \Delta^{1/2}$ is the polar decomposition 
of $S$, then $\Delta$ is a positive selfadjoint, $J$ is an anti-unitary involution 
and $J \Delta J = \Delta^{-1}$. 
\item[\rm(iii)] $S^* = J S J$ is the complex conjugation with respect to the 
standard subspace $J(V)$. 
\item[\rm(iv)] The real orthogonal complement of $V$ w.r.t.\ $\Re \la \cdot,\cdot\ra$ is 
$iJ(V)$. In particular, $\cH = V \oplus i J(V)$ as a real Hilbert space. 
\item[\rm(v)] $J(V) = V^{\bot,\omega}$ for the symplectic form $\omega := \Im \la \cdot,\cdot\ra$. 
\end{itemize}
\end{lem}

\begin{prf} (i), (ii):  
For $x,y\in V$ and $z := x + iy$, we have 
\begin{equation}
  \label{eq:gnorm}
 \|(z,S(z))\|^2 
= \|x + i y\|^2 + \|x - i y\|^2 
= 2(\|x\|^2 + \|y\|^2),
\end{equation}
so that the graph of $S$ is closed because $V \oplus V$ is complete, which in turn follows from 
the closedness of $V$. 
We conclude that $S$ is a closed densely defined operator on $\cH$, 
so that $S$ has a polar decomposition as in (ii). 
Since $\im(S) = V_\C$ is dense, $J$ is an isometry. 
With the same arguments as in \cite[Prop.~2.5.11]{BR02}, it now follows that 
$J$ is an anti-unitary involution satisfying 
$J \Delta  J = \Delta^{-1}.$

(iii) This further leads to 
$S^* = \Delta^{1/2} J  = J S J,$ 
which shows that $S^* = S_{J V}$ is the complex conjugation corresponding to the 
standard subspace $J(V)$. 

(iv) For $w \in \cH$, we have the following chain of equivalences: 
\begin{align*}
 \Re\la iV, w \ra = \{0\} 
\Leftrightarrow & (\forall v \in V_\C) \Re \la v, w \ra = \Re \la S v, w \ra = \Re \la w, S v\ra\\\Leftrightarrow & (\forall v \in V_\C) \la v, w \ra =\la w, S v\ra \ 
\Leftrightarrow\  w \in \cD(S^*) \ \& \ S^*w = w \ \Leftrightarrow \ w \in J(V). 
\end{align*}
(v) follows immediately from (iv) because 
$i V^{\bot,\omega}$ is the orthogonal complement w.r.t.\ $\Re \la \cdot,\cdot\ra$. 
\end{prf}

\begin{rem} By Lemma~\ref{lem:e.1}, the scalar product $\la \cdot,\cdot \ra$ 
is real-valued on $V \times V$ if and only if $V \subeq J(V)$, which is equivalent to $J(V) = V$. 
This happens only if $S^* = S$, and this in turn is equivalent to $\Delta = \1$, i.e., $J = S$. 
\end{rem}

\begin{lem} \mlabel{lem:e.2} Let $\Delta$ be an injective 
 positive selfadjoint operator on 
$\cH$ and $J$ an anti-unitary involution with 
$J \Delta J = \Delta^{-1}$. Then 
\begin{itemize}
\item[\rm(i)] $S := J \Delta^{1/2} \: \cD(\Delta^{1/2}) \to \cH$ is a closed antilinear 
involution and $V := \{ v \in \cD(S) \: Sv = v\}$ is a standard real subspace of 
$\cH$ with $S_V = S$. 
\item[\rm(ii)] The selfadjoint operator $H := \log\Delta$ satisfies $J H J= -H$ 
and $U_{(t,\eps)} := e^{itH} J^\eps = \Delta^{it} J^\eps$ 
defines an anti-unitary representation of the direct product 
group $\R \times \Z/2\Z$. 
\end{itemize}
\end{lem}

\begin{prf} That $S$ is an involution follows from 
$S^2 = J\Delta^{1/2} J \Delta^{1/2} = \Delta^{-1/2} \Delta^{1/2} = \id_{\cD(S)}$. 
Further, the closedness of the selfadjoint operator $\Delta^{1/2}$ implies 
that $S$ is closed. Now (i) follows from \eqref{eq:gnorm}, and (ii) is clear. 
\end{prf}

\begin{ex} \mlabel{ex:4.5} (Modular data of von Neumann algebras) 
Tomita--Takesaki theory (\cite[\S 2.5]{BR02}) 
starts with a cyclic separating vector $\Omega$ 
for a von Neumann algebra $\cM \subeq B(\cH)$. 
Then $V := \oline{\{M\Omega \: M^* = M \in \cM\}}$ is a standard real subspace of $\cH$, 
and for $S = S_V$, $\Delta := S^*S$ is the {\it modular  operator} and 
$\alpha_t(M) := \Delta^{-it} M \Delta^{it}$ 
defines a one-parameter group of automorphisms of $\cM$ (the {\it modular group}). 
Further, the {\it modular involution} $J$ satisfies $J \cM J = \cM'$. 
\end{ex}

The preceding two lemmas imply that we have a one-to-one correspondence between 
the following data. 
\begin{itemize}
\item[\rm(a)] standard closed subspaces $V \subeq \cH$. 
\item[\rm(b)] pairs $(\Delta, J)$, where $0 < \Delta = \Delta^*$ and $J$ is an 
anti-unitary involution satisfying $J\Delta J = \Delta^{-1}$. 
\item[\rm(c)] pairs $(H, J)$, where $H = H^*$ and $J$ is an 
anti-unitary involution satisfying $JHJ = -H$. 
\item[\rm(d)] strongly continuous anti-unitary representations of $\R \times \Z/2\Z 
\cong \R^\times$. 
\end{itemize}

\begin{rem}  \mlabel{rem:4.6} 
Let $U^c_t = e^{itH}$  be a unitary one-parameter group on $\cH$ and  
$J$ be an anti-unitary involution on $\cH$ with $J H J = -H$. 
Then $S := J e^{-\beta H/2}$ is an unbounded antilinear involution with 
$\cD(S) = \cD(e^{-\beta H/2})$. We further assume that $\cV$ is a real vector space and that 
$j \: \cV \to \cH$ is a linear map with cyclic range contained in 
$\Fix(S) = \{ v \in \cD(S) \: Sv =v\}$. As the operators $U^c_t$ commute with 
$S$, $\Fix(S)$ is a closed $U^c$-invariant subspace. 
Now $\phi(z)(v,w) := \la  e^{-z H} j(v), j(w)\ra$ defines for $0 \leq \Re z \leq \beta$ 
an element $\phi(z) \in \Sesq(\cV)$. 
It is easy to see that the kernel 
$\phi\big(\frac{z+\oline w}{2}\big)$ is positive 
definite. For $v,w \in \cV$, we further have: 
\begin{eqnarray}\label{eq:phi-rel}
\phi(z)(v,w) 
&=& \la e^{-zH} j(v), j(w) \ra 
= \la e^{-zH} Sj(v), Sj(w) \ra 
= \la e^{-zH} J e^{-\beta H/2}j(v), J e^{-\beta H/2}j(w) \ra \notag\\
&=& \la J e^{\oline zH} e^{-\beta H/2}j(v), J e^{-\beta H/2}j(w) \ra 
= \la  e^{-\beta H/2}j(w) , e^{\oline zH} e^{-\beta H/2}j(v) \ra \notag \\
&=& \la  j(w), e^{-(\beta-\oline z) H}j(v) \ra 
= \oline{\phi(\beta-\oline z)(v,w)}.
\end{eqnarray}
This means that, for $0 \leq \Re z \leq \beta$, the relation 
 $\phi(\beta-\oline z) = \oline{\phi(z)}$ holds in $\Sesq(\cV)$. 
From that we derive in particular that $\Re \phi$ defines a reflection 
positive function $\T_\beta \to \Herm(\cV)$. Note that all functions 
$\phi^{v,v}(t) := \phi(t)(v,v)$ are real-valued, hence reflection positive. 
The function $\phi$ itself admits 
a $\beta$-periodic extension to $\R$ if and only if $\phi(0)$ is real-valued on 
$V \times V$. Since the scalar product of $\cH$ is in general not real-valued on $V \times V$, 
the form $\phi(0)$ need not be symmetric 
\end{rem}

\section{KMS states} 
\mlabel{sec:5}

Let $(\cA, \R, \alpha)$ be a $C^*$-dynamical system, i.e., 
$\alpha \: \R \to \Aut(\cA)$ is a homomorphism defining a continuous $\R$-action on 
the $C^*$-algebra $\cA$. 
We recall from \cite[Props.~5.3.3, 5.3.7]{BR96} that an $\alpha$-invariant 
state $\omega$ of $\cA$ is a {\it KMS state at value $\beta > 0$} if, 
for all $A,B\in \cA$, there exists a bounded holomorphic function 
$F_{A,B} : \cD_\beta \to \C$ 
extending continuously to $\oline{\cD_\beta}$, such that 
\begin{equation}
  \label{eq:kms-def}
F_{A,B}(t) = \omega(A \alpha_t(B)) \quad \mbox{ and } \quad 
F_{A,B}(t + i \beta) = \omega(\alpha_t(B)A) \quad \mbox{ for  } \quad t \in \R.
\end{equation}

The following lemma links our concept of a positive definite function 
with values in sesquilinear forms to invariant states of $\cA$ (Definition~\ref{def:a.3}). 
\begin{lem} \mlabel{lem:5.1} {\rm(a)} Let $(\cA,G,\alpha)$ be a $C^*$-dynamical system and 
$\omega$ be a $G$-invariant state of $\cA$. Then 
\[ \psi \: G \to \Sesq(\cA), \quad 
\psi(g)(A,B) := \omega(B^*\alpha_g(A)) \] 
is a positive definite function.
\end{lem}

\begin{prf} Let $(\pi_\omega, U, \cH_\omega, \Omega)$
be the cyclic covariant representation with $U_g \Omega= \Omega$ for $g \in G$ 
(\cite[Cor.~2.3.17]{BR02}). For the $\R$-equivariant linear map 
$j \: \cA \to \cH_\omega, j(A) = \pi(A)\Omega$ with dense range we then have 
\[ \psi(g)(A,B) 
= \la \pi(B^*\alpha_g(A))\Omega, \Omega \ra 
= \la U_g \pi(A) \Omega, \pi(B)\Omega \ra
= \la U_g j(A), j(B)\ra,\] 
Therefore the assertion follows from Proposition~\ref{prop:gns}.
\end{prf}

\begin{prop}
For $G = \R$, an invariant state $\omega$ is a KMS state if and only if 
$\psi$ extends to a function $\oline{\cD_\beta} \to \Sesq(\cA)$ which is pointwise 
continuous and pointwise holomorphic on $\cD_\beta$ on $\cA \times \cA$ and satisfies 
\begin{equation}
  \label{eq:kms-ancont}
\psi(t + i \beta)(A,B) = \psi(-t)(B^*,A^*) = \oline{\psi(t)(A^*,B^*)} 
\quad \mbox{ for } \quad t \in \R.
\end{equation}
\end{prop}

\begin{prf} The fact that $\psi$ is positive definite 
implies in particular that $\psi(-t) = \psi(t)^*$: 
\[ \psi(-t)(A,B) = \omega(B^* \alpha_{-t}(A)) = \omega(\alpha_t(B^*)A) 
= \oline{\omega(A^* \alpha_t(B))} = \psi(t)^*(A,B).\] 
Comparing with \eqref{eq:kms-def}, we see that $\psi(t)(A,B) = F_{B^*,A}(t)$, so that the 
assertion follows from $\omega(\alpha_t(A)B^*) = \omega(A \alpha_{-t}(B^*)) = \psi(-t)(B^*,A^*)$ 
and \cite[Props.~5.3.3, 5.3.7]{BR96}. 
\end{prf}

\begin{rem} \mlabel{rem:5.3} 
(a) By analytic continuation, the relation \eqref{eq:kms-ancont} implies that 
\[ \psi(i \beta-z)(A,B) = \psi(z)(B^*,A^*) = \oline{\psi(-\oline z)(A^*,B^*)} 
\quad \mbox{ for } \quad z \in \oline{\cD_\beta}.\] 
For $\phi(t) := \psi(it)$, we obtain in particular 
\[ \phi(\beta-t)(A,B) = \phi(t)(B^*,A^*) = \oline{\phi(t)(A^*,B^*)} 
\quad \mbox{ for } \quad 0 \leq t \leq \beta.\] 
Using the notation from \ref{sec:a}, we define 
$\phi^{A,B}(t) := \phi(t)(A,B)$. Then 
$\phi^{A,A}(t) := \phi(t)(A,A)$ is real for $A \in \cA$ and that 
$\phi^{A,A^*}(\beta - t) = \phi^{A,A^*}(t)$. For $A = A^*$, it follows in particular 
that the function $\phi^{A,A}$ on $[0,\beta]$ defines a reflection positive function 
on $\T_\beta$. 

(b) In the situation of the proof of Lemma~\ref{lem:5.1}(b), for $G = \R$,  
Lemma~\ref{lem:a.2} implies that $\pi(\cA)\Omega \subeq \cD(e^{-\beta H/2})$ and that 
\[ \phi(t)(A,B) = \la e^{-tH} \pi(A)\Omega, \pi(B)\Omega \ra \quad \mbox{ for } \quad 
0 \leq t \leq \beta.\] 

(c) In the context of Tomita--Takesaki theory (Example~\ref{ex:4.5}), 
we put $U_t := \Delta^{-it}$ and $H = - \log \Delta$, so that 
$\Delta^{1/2} = e^{-H/2}$, which corresponds to $\beta = 1$. 
For the state $\omega(M) := \la M\Omega, \Omega\ra$ 
and $\psi(t)(A\Omega, B\Omega) = \la U_t A\Omega, B \Omega \ra = \phi(-it)(A\Omega, B\Omega)$, 
we obtain from \eqref{eq:phi-rel} for $A = A^*, B = B^*$
\[ \psi(t+ i)(A\Omega,B\Omega) 
= \phi(1- it)(A\Omega,B\Omega)  
= \oline{\phi(-it)}(A\Omega,B\Omega)  
= \oline{\psi(t)}(A\Omega,B\Omega),\] 
which is the KMS condition \eqref{eq:kms-ancont} for $\beta = 1$ 
(\cite[Prop.~5.3.10]{BR96}).
\end{rem}

\section{Perspectives} 
\mlabel{sec:6}

The main outcome of the present note is that it clarifies the connection between 
reflection positivity for the triple $(\T_\beta, \T_\beta^+, \tau_\beta)$ and the modular data 
$(J,\Delta)$ arising naturally in Tomita--Takesaki theory, resp., the theory of KMS states. 
We hope that this will serve as a basis for a deeper understanding of 
reflection positivity for higher dimensional groups. 

{\bf Relativistic KMS states:} One interesting direction is to connect the 
relativistic KMS states for a $d$-dimensional space time 
introduced by J.~Bros and D. Buchholz in \cite{BB94} to 
reflection positivity for the group $\R^d$ 
(see also \cite{GJ06} for $d = 2$). 
Here the strip $\cD_\beta \subeq \C$ is replaced 
by tube domain 
\[ \cT_\beta := \{ z \in \C^d \: \Im z \in V_+ \cap (\beta e - V_+) \}, \] 
where $V_+$ is the open forward light cone and $e \in V_+$ is a timelike vector of unit length. 
One expects a duality between $G^c = \R^d$ and $G = \T_\beta \times \R^{d-1}$ and a suitable 
domain $G^+\subeq G$ for reflection positivity. 

{\bf ax+b-group and generalizations:} There are still several aspects of reflection positivity for 
the $ax+b$-group $G \cong \R \rtimes \R^\times$ (the affine group of the real line) that 
are not covered by the discussion in \cite{NO14b}. 
Here the only non-trivial involution is given by 
$\tau(b,a) = (-b,a)$. In \cite{NO14b} we have seen how reflection positivity 
for $(G,G^+, \tau)$, where $G^+ = \{ (b,a) \: b > 0, a >0\}$, leads to euclidean realizations  
of those representations $\pi^c$ of $G^c \cong \R \rtimes \R^\times_+$ satisfying a 
positive spectrum condition for translations. For any anti-unitary involution 
$J$ satisfying $J \pi(b,a) J = \pi(-b,a)$, these representations extend 
naturally by $\pi(0,-1) := J$ to anti-unitary 
representations of the non-connected group $\R \rtimes \R^\times$.  
Such representations arise naturally in the context of Borchers' triples 
(\cite{BLS11}, \cite{Bo92}), where they actually extend to anti-unitary representations 
of $G^c := \R^d \rtimes_\gamma \R^\times \cong (\R^2 \times \R^{d-2}) \rtimes \SO_{1,1}(\R)$ 
and $J = \pi(0,0,-1)$ acts on Minkowski space $\R^d \cong \R^2 \oplus \R^{d-2}$ 
by $J \pi(b,c,\1)J  = \pi(-b,c,\1)$, where $\R^2$ is $2$-dimensional Minkowski space and 
$\R^{d-2}$ is space-like. The dual group is the subgroup 
$G \cong (\R^2 \times \R^{d-2}) \rtimes \SO_2(\R)$ of the $d$-dimensional motion group. 
This should provide a natural framework for extending the reflection 
positivity for the circle discussed here in which one can also treat relativistic 
KMS states. Clearly, the crucial case to be understood is $d= 2$. 
On $G\cong \R^d \rtimes \SO_2(\R)$, 
the involution induced by the time reflection $\theta$ is given by 
$\tau(b,a) = (\theta(b),a^{-1})$ and likewise on $G^c$. Note that this involution 
commutes with the action of $J$. 

{\bf Reflection positivity on the side of $G^c$:} Classically, reflection positivity 
is a condition imposed in the euclidean model where we have a unitary representation 
$(\pi,\cE)$ of the group $G$, whereas on the dual side we have a representation 
$(\pi^c, \hat\cE)$ which has no specific reflection symmetry. However, 
in the theory of modular inclusions, one encounters anti-unitary representations of the 
$ax+b$-group corresponding to a modular pair $(J,\Delta)$ and a unitary one-parameter group 
$U$ which satisfy the commutation relations 
\[ \Delta^{it} U(s) \Delta^{-it} = U(e^{-2\pi t}s) \quad \mbox{ and } \quad 
J U(s) J = U(-s) \quad \mbox{ for } \quad t,s \in \R.\] 
In \cite[Thm.~3.2]{BGL02} it is shown that, in this context, the infinitesimal 
generator $H$ of $U(s) = e^{isH}$ satisfies $H \geq 0$ if and only if 
$U(s)V \subeq V$ holds for $s \geq 0$, where $V$ is 
the standard real subspace specified by 
$(J,\Delta)$ (cf.\ Lemma~\ref{lem:e.2}). This situation has a remarkable similarity 
with the reflection positivity for $(\R,\R_+,-\id_\R)$, where we have a unitary 
one-parameter group $(U_t)_{t \in \R}$ of $\cE$ whose positive part acts by 
isometries on the subspace $\cE_+$ and the spectrum of the generator of the dual 
one-parameter group $U^c$ is positive (\cite{NO14b}). However, here the subspace 
$V$ is real, $\la Jv,v\ra \geq 0$ for $v \in V$ and $(U(s))_{s \geq 0}$ acts by 
real isometries on $V$. 

{\bf More general reflection positive functions:} For a symmetric 
Lie group $(G,\tau)$, the notion of a reflection positive function $\phi \: G \to \C$ 
has been introduced with the idea that these should be positive definite functions 
corresponding to unitary representations $U$ on a reflection positive Hilbert space 
$(\cE,\cE_+, \theta)$ such that $\phi(g) = \la \pi(g)v,v\ra$ holds for a 
$\theta$-fixed vector $v \in \cE_+$. Then reflection positivity with respect to a domain 
$G^+ \subeq G$ corresponds to $\pi(G^+)^{-1}v \subeq \cE_+$ (cf.\ \cite{NO14a}) 
and is encoded in the positive definiteness of the kernel $\phi(g\tau(h)^{-1})$ on 
$G^+$. For $(\R,\R_+, -\id_\R)$ and $(\T_\beta,\T_\beta^+,\tau_\beta)$, 
this implies that 
$\phi$ is real-valued (cf.\ Theorem~\ref{thm:2.4} for $\dim\cV =1$). 
But this condition is too restrictive to cover the 
$\beta$-periodic functions $\phi^{A,A^*}$ 
from Remarks~\ref{rem:4.6} (see also Remark~\ref{rem:5.3}). 
This suggests to work with a more general concept where we give up the 
$\theta$-invariance of $v$, so that $\phi$ need no longer be $\tau$-invariant. 
Extending $\pi \: G \to \U(\cE)$ to a representation $\hat\pi$ 
of the group $G_\tau := 
G \rtimes \{\1,\tau\}$ by $\pi(\1,\tau) := \theta$, the condition 
$\pi(G^+)^{-1} v \subeq \cE_+$ corresponds to the positive definiteness of the kernel 
$\la \theta\pi(h^{-1})v,\pi(g^{-1})v \ra = \la \hat\pi(g\tau(h)^{-1},\tau)v,v\ra$, 
which can be expressed in terms of the positive definite 
function $\hat\phi(g) := \la \hat\pi(g)v,v\ra$ on the non-connected group $G_\tau$. 
It is therefore desirable to obtain an explicit description of the 
so specified functions 
for $(\R,\R_+, -\id_\R)$ and $(\T_\beta, \T_\beta^+, \tau_\beta)$ and thus 
a generalization of the corresponding results in \cite{NO14a} for $\R$ and 
Theorem~\ref{thm:2.4} for $\T_\beta$. 
Of course, it would be of particular interest to see if 
the functions 
$\phi^{A,A^*}$ from Remarks~\ref{rem:4.6} fall into this larger class. 

\ack{K.-H. Neeb was supported by DFG-grant NE 413/7-2, Schwerpunktprogramm 
``Darstellungstheorie'' and Gestur \'Olafsson  was supported by NSF grant 
DMS-1101337 and the Emerging Fields Project 
``Quantum Geometry'' of the University of Erlangen.}

\appendix
\def\thesection{Appendix \Alph{section}}   

\section{Form-valued kernels} 
\mlabel{sec:a}

\def\thesection{\Alph{section}}   

Classically, reproducing kernels arise from Hilbert spaces
$\cH$ of functions $f \: X \to \C$ ($X$ a set) 
for which the evaluations $f \mapsto f(x)$ are continuous, 
hence representable by elements $K_x \in \cH$ by $f(x) = \la f, K_x\ra$, and 
then $K(x,y) := K_y(x) = \la K_y, K_x\ra$ is called the {\it reproducing kernel} of the space 
$\cH$. Then $K$ determines $\cH$ uniquely. Accordingly, we  write 
$\cH_K\subeq \C^X$ for the Hilbert space determined by $K$. It is a classical 
result that a kernel function $K \: X \times X \to \C$ is the reproducing kernel 
of some Hilbert space if and only if it is {\it positive definite} in the sense 
that, for any finite collection $x_1, \ldots, x_n \in X$, the matrix 
$(K(x_j,x_k))_{1 \leq j,k\leq n}$ is positive semidefinite (cf.\ \cite{Ar50}, \cite[Ch.~1]{Ne00}). 
There is a natural generalization to Hilbert spaces $\cH$ of functions 
with values in a Hilbert space $\cV$, i.e., $\cH \subeq \cV^X$. Then $K_x(f) = f(x)$ is a 
linear operator $K_x \: \cH \to \cV$ and we obtain a kernel 
$K(x,y) := K_x K_y^* \in B(\cV)$ with values in the bounded operators on $\cV$. 
However, there are also situations where one would like to deal with kernels 
whose values are unbounded operators, so that one has to generalize this 
context further. As we shall explain below, the concept of a positive 
definite kernel with values in the space $\Sesq(\cV)$ of sesquilinear complex-valued 
forms on $\cV$ provides a natural context to deal with all relevant cases. 
\begin{footnote}{If $\cV$ is a real space, then $\Sesq(\cV)$ is the 
space of bilinear maps $\cV \times \cV \to \C$, and if $\cV$ is complex, 
it stands for those maps $\cV \times \cV \to \C$ 
which are linear in the first and antilinear in the second 
argument. The uniqueness of sesquilinear extension leads to 
$\Sesq(\cV) = \Sesq(\cV_\C)$ for a real space $\cV$. If $\cV$ is complex, then polarization shows that 
every $\phi \in \Sesq(\cV)$ is uniquely determined by its values on the diagonal.} 
  \end{footnote}

For a real or complex vector space $\cV$, we write $\cV^\sharp$ for the complex vector space of 
antilinear maps $\cV \to \C$. 

\begin{defn} Let $X$ be a set. We call a map 
$K \: X \times X \to \Sesq(\cV)$ 
{\it a positive definite kernel} if the associated scalar-valued kernel 
\[  K^\flat \: (X \times \cV) \times (X \times \cV) \to \C, \quad 
K^\flat((x,v),(y,w)) := K(x,y)(w,v) \] 
is positive definite. 
For elements $f$ of the corresponding 
reproducing kernel Hilbert space $\cH_{K^\flat} \subeq \C^{X \times \cV}$,  
we then have 
\[ f(y,w) = \la f, K^\flat_{y,w}\ra \quad \mbox{ with } \quad 
K^\flat_{y,w}(x,v) := K^\flat((x,v),(y,w)) = K(x,y)(w,v).\] 
Therefore $v \mapsto K_{x,v}^\flat$ is linear, and this implies that $f(x,\cdot)$ is antilinear. 
We identify $\cH_{K^\flat}$ with a subspace of $(\cV^\sharp)^X$ by identifying 
$f \in \cH_{K^\flat}$ with the function $f^\sharp \: X \to \cV^\sharp, f^\sharp(x):= f(x,\cdot)$. 
We call 
\[ \cH_K := \{ f^\sharp \: f \in \cH_{K^\flat} \} \subeq (\cV^\sharp)^X \] 
the {\it (vector-valued) reproducing kernel space associated to $K$}. The elements 
\[ K_{x,v} := (K^\flat_{x,v})^\sharp =K(\cdot, x)(v,\cdot), \quad x \in  X, v \in \cV, \] 
then form a dense subspace of $\cH_K$ with 
\begin{equation}
  \label{eq:app}
\la K_{y,w}, K_{x,v} \ra =K^\flat((x,v),(y,w)) = K(x,y)(w,v).
\end{equation}
\end{defn}

\begin{ex} \mlabel{ex:a.2} (a) If $K \: X \times X \to B(\cV)$ is an operator-valued kernel, 
where $\cV$ is a complex Hilbert space, 
then we obtain a $\Sesq(\cV)$-valued kernel by $Q(x,y)(v,w) := \la K(x,y)v, w \ra$. 
We have a natural inclusion $B(\cV)\into \Sesq(\cV), A \mapsto \la A\cdot, \cdot\ra$, 
whose range is the space of continuous sesquilinear forms. All the functions 
$f^\sharp \: X \to \cV^\sharp$ in $\cH_Q$ take values in continuous functionals, hence can be identified 
with $\cV$-valued functions. This leads to the realization $\cH_K \into \cV^X$ 
(cf.\ \cite[Ch.~1]{Ne00}). 

(b) For a one-point set $X = \{*\}$, a positive definite $\Sesq(\cV)$-valued kernel 
is simply a positive semidefinite $K \in \Sesq(\cV)$. 
The corresponding Hilbert space $\cH_K \subeq \cV^\sharp$ is generated by the elements 
$K_v = K(v,\cdot)$ with $\la K_v, K_w \ra = K(v,w)$. 

In particular, if $\cV$ is a Hilbert space, then the natural inclusion 
$\cV \into \cV^\sharp, v \mapsto \la v, \cdot\ra$, corresponds to the kernel 
$K(v,w) = \la v, w\ra$. 

(c) If $\cV = \C$, then $\Sesq(\cV) \cong \C$ and $\Sesq(\cV)$-valued kernels are 
complex-valued kernels. 

(d) If $\cA$ is a $C^*$-algebra and $\omega \in \cA^*$ a positive functional, then 
$K_\omega(A,B) :=\omega(B^*A)$ 
is a positive semidefinite sesquilinear kernel for which the corresponding 
Hilbert space $\cH_{K_\omega}\subeq \cA^\sharp$ can be obtained from the 
GNS representation $(\pi_\omega, \cH_\omega, \Omega)$ (\cite[Cor.~2.3.17]{BR02}) by the embedding 
\[ \Gamma \: \cH_\omega \to \cA^\sharp, \quad 
\Gamma(\xi)(A) := \la \xi, \pi(A)\Omega \ra \quad \mbox{ because } \quad 
\la \pi(A)\Omega, \pi(B)\Omega \ra = \omega(B^*A) = K_\omega(A,B).\] 
Note that $\cA$ has a natural representation on $\cA^\sharp$ by 
$(A.\beta)(B) := \beta(A^*B)$ and that $\Gamma$ is equivariant with respect to this 
representation.\begin{footnote}{This realization of the Hilbert space 
$\cH_\omega$ has the advantage that we can see its elements as elements of the 
space $\cA^\sharp$ (see \cite{Ne00} for many applications of this perspective). 
Usually, $\cH_\omega$ is obtained as the Hilbert completion of 
a quotient of $\cA$ by a left ideal which leads to a much less concrete space.   
}\end{footnote}
\end{ex}

\begin{defn} \mlabel{def:a.3} Let $G$ be a group. 
A function $\phi \: G \to \Sesq(\cV)$ is said to be {\it positive definite} 
if the $\Sesq(\cV)$-valued kernel $K(g,h) := \phi(gh^{-1})$ is positive definite. 
\end{defn}

\begin{prop}{\rm(GNS-construction)} \mlabel{prop:gns} 
{\rm(a)} Let $\phi \: G \to \Sesq(\cV)$ be a positive definite function. Then 
$(\pi_\phi(g)f)(h) := f(hg)$ defines a unitary representation of $G$ on the 
reproducing kernel Hilbert space $\cH_\phi \subeq (\cV^\sharp)^G$ with 
kernel $\phi(gh^{-1})$ and the range of the map $j \: \cV \to \cH_\phi, j(v)(g):=  
\phi(g)(v,\cdot)$ is a cyclic subspace, i.e., $\pi(G)j(\cV)$ spans a dense subspace of $\cH$. 
We then have 
\[ \phi(g)(v,w) = \la \pi_\phi(g)j(v), j(w) \ra \quad \mbox{ for } \quad g \in G, v,w, \in \cV.\] 

{\rm(b)} If, conversely, $(\pi, \cH)$ is a unitary representation of $G$ and 
$j \: \cV \to \cH$ a linear map whose range is cyclic, then 
\[ \phi \: G \to \Sesq(\cV), \quad \phi(g)(v,w) := \la \pi(g)j(v), j(w) \ra \] 
is a $\Sesq(\cV)$-valued positive definite function and 
$(\pi, \cH)$ is unitarily equivalent to $(\pi_\phi, \cH_\phi)$. 
\end{prop}

\begin{prf} (cf.\ \cite[Sect.~3.1]{Ne00}) (a) For the kernel 
$K(g,h) := \phi(gh^{-1})$ and $v \in \cV$, the right invariance of the kernel $K$ implies that 
$\pi(g)K_{h,v} = K_{hg^{-1}, v}$, and 
from that one easily derive the invariance of $\cH_\phi$ under right translations 
and the unitarity of their restrictions to $\cH_\phi$. 
Then $j(v) = \phi(\cdot)(v,\cdot) = K_{\1,v}$ and 
$\pi(g)j(v) = K_{g^{-1},v}$ imply that $j(\cV) \subeq \cH_\phi$ is cyclic. 
Finally we note that 
\[ \la \pi(g)j(v), j(w) \ra = \la K_{g^{-1},v}, K_{\1,w} \ra = K(g)(v,w) = \phi(g)(v,w).\]

(b) The positive definiteness of $\phi$ follows easily from the relation 
$\phi(gh^{-1})(v,w) = \la \pi(h)^{-1}j(v), \pi(g)^{-1}j(w)\ra$. Since $j(\cV)$ is cyclic, 
the map 
\[ \Gamma(\xi)(g)(v) := \la \xi, \pi(g)^{-1} j(v) \ra \]
defines an injection $\cH \into (\cV^\sharp)^G$ whose range is the subspace 
$\cH_\phi$ and which is equivariant with respect to the right translation action 
$\pi_\phi$. 
\end{prf}

\begin{rem} If $\cV$ is a Hilbert space 
and $j$ is continuous, then we have the adjoint operator 
$j^* \: \cH \to \cV$ is well-defined and we obtain 
the $B(\cV)$-valued positive definite function 
$\phi(g) := j^* \pi(g) j$ which can be used to realize $\cH$ in $\cV^G$ 
(cf.\ Example~\ref{ex:a.2}(a)). 
\end{rem}

\def\thesection{Appendix \Alph{section}}   
\section{Integral representations} 
\def\thesection{\Alph{section}}   
\mlabel{sec:b}

For a more concrete realization of unitary representations associated to 
positive definite functions in $L^2$-spaces, integral representations are of crucial importance. 
The following result is a straight-forward generalization of Bochner's Theorem 
for locally compact abelian groups. Here we write $\Sesq(\cV)_+ \subeq \Sesq(\cV)$ for the 
convex cone of positive semidefinite forms. 

\begin{prop} \mlabel{prop:bochner} Let $G$ be a locally compact abelian group. 
If $\phi \: G \to \Sesq(\cV)$ is a positive definite function 
for which all functions $\phi^{v,w} := \phi(\cdot)(v,w), v,w \in \cV$,
 are continuous, then there exists a uniquely determined $\Sesq(\cV)_+$-valued Borel 
measure $\mu$ on the locally compact group $\hat G$ 
such that $\hat\mu(g) := \int_{\hat G} \chi(g)\, d\mu(\chi) = \phi(g)$ holds 
for every $g \in G$ pointwise on $\cV \times \cV$. 
\end{prop}

\begin{prf} First, Bochner's Theorem for scalar-valued positive definite 
functions yields for every $v \in \cV$ a finite positive measure $\mu^v$ on $\hat G$ such that 
$\phi^{v,v} = \hat{\mu^v}$. By polarization, we obtain complex measures 
$\mu^{v,w}$, $v,w \in \cV$, on $\hat G$ with $\phi^{v,w} = \hat{\mu^{v,w}}$. Then the 
collection $(\mu^{v,w})_{v,w \in \cV}$ of complex measures on $\hat G$ defines a 
$\Sesq(\cV)_+$-valued measure by 
$\mu(\cdot)(v,w) := \mu^{v,w}$ for $v,w \in \cV,$ 
and this measure satisfies $\hat\mu = \phi$. 
\end{prf} 

\begin{rem} Suppose that $E$ is the spectral measure on the character group $\hat G$ for which 
the continuous unitary representation $(\pi, \cH)$ is represented by 
$\pi(g) = \int_{\hat G} \chi(g)\, dE(\chi)$. Then, for $v \in \cH$, the positive definite 
function $\pi^v(g) := \la \pi(g)v,v\ra$ is the Fourier transform of the measure 
$E^{v,v} = \la E(\cdot)v,v\ra$. This establishes a close link between spectral measures 
and the representing measures in the preceding proposition.   
\end{rem}

We say that a subset $\cD$ in the real vector space $E$ (contained in 
$E_\C = E + i E$) is {\it
finitely open} if, for each finite-dimensional subspace $F \subeq E$,
the intersection $F \cap \cD$ is open in $F$. If $\cV$ is a Hilbert
space, then a function 
$f \: T_\cD := \cD + i E\to E$ on the corresponding tube domain is called {\it
Gateaux holomorphic} if the restriction to each tube domain 
$T_{F \cap \cD}$ is holomorphic. We write 
$\cO_G(T_\cD,\cV)$ for the space of all Gateaux holomorphic
$\cV$-valued functions on $T_\cD$. For the theory of holomorphic
functions on domains in infinite dimensional spaces we refer to Herv\'e's monograph 
\cite{He89} or \cite{Ne00} for the connections to representation theory. 

\begin{thm}\mlabel{thm:I.7} {\rm(Laplace transforms and positive definite kernels)} 
Let $E$ be a real vector space and $\cD \subeq E$ 
a non-empty convex finitely open subset. Let $\cV$ be a Hilbert space and 
$\phi \: \cD \to B(\cV)$ such that 
\begin{itemize}
\item[\rm(L1)] the kernel $Q_\phi(x,y) 
=\phi\big(\frac{x+y}{2}\big)$ is positive definite. 
\item[\rm(L2)] $\phi$ is weak operator continuous on very line segment in $\cD$, 
i.e., all functions \break $t \mapsto \la \phi(x+th)v,v\ra$, $v \in \cV$, are continuous 
on $\{ t \in \R \: x + th \in \cD\}$. 
\end{itemize}
Then the following assertions hold: 
\begin{itemize}
\item[\rm(i)] There exists a unique $\Herm^+(\cV)$-valued measure 
$\mu$ on the smallest $\sigma$-algebra on $E^*$ for which all point
evaluations are continuous such that 
\[ \phi(x) = \cL(\mu)(x) := \int_{E^*} e^{-\lambda(x)}\, d\mu(\lambda) \quad \mbox{ for } \quad 
x \in \cD.\] 
\item[\rm(ii)] The map 
\[ {\cal F} \: L^2(E^*, \mu;\cV) \to \cO_G(T_\cD,\cV),\quad 
\la \cF(f)(z), v \ra  := \la f, e_{-\frac{\oline z}{2}}v \ra \] 
is unitary  onto the reproducing kernel space ${\cal H}_\phi := \cH_{Q_\phi}$
corresponding to the kernel associated to $\phi$. 
It intertwines the unitary representation 
\[ (\pi(u)f)(\alpha) := e^{i\alpha(u)}f(\alpha) \quad \mbox{ on } \quad 
L^2(E^*,\mu) \quad \mbox{ and } \quad 
(\tilde\pi(u)f)(z) := f(z - 2 i u)\quad \mbox{ on } \quad \cH_\phi. \] 
\item[\rm(iii)] $\phi$ extends to a  unique 
Gateaux holomorphic function $\hat\phi$ on the tube domain $T_\cD$ 
which is positive definite in the sense that the
 kernel $\hat\phi\big(\frac{z + \oline w}{2}\big)$ is positive definite.
\end{itemize}
\end{thm}

\begin{prf} (i) follows from \cite[Thm.~18.8]{Gl03} and (ii) from \cite[Thm.~I.7]{NOr02}, provided 
we show that (L2) implies that, for every trace class operator $S \in B_1(\cV)$, the function 
$x \mapsto \tr(\phi(x)S)$ is continuous on line segments in $\cD$. 
To this end, we may assume that $E = \R$ and $\cD = ]a,b[$. 
Since $B_1(\cH)$ is spanned by positive operators with trace $1$, we may w.l.o.g.\ 
assume that $S = S^* \geq 0$ with $\tr S = 1$. We may further assume that 
$\dim \cV = \infty$; otherwise the assertion is trivial. 
Then there exists an orthogonal sequence $(v_n)_{n \in \N}$ in $\cV$ with 
$Sv = \sum_n \la v, v_n \ra v_n.$ 
This leads to 
\[ \tr(\phi(x)S) = \sum_n \phi_n(x) \quad \mbox{ with } \quad 
\phi_n(x) := \la \phi(x)v_n, v_n \ra.\] 
By assumption, all functions $\phi_n$ are continuous. 
Applying \cite[Thm.~I.7]{NOr02} to the case $\cV = \C$, we obtain 
positive measures $\nu_n$ on $\R$ with with 
$\phi_n = \cL(\nu_n)$ on $]a,b[$. Let $\nu := \sum_n \nu_n$. Then 
\[ \cL(\nu)(x) = \sum_n \cL(\nu_n)(x) = \sum_n \phi_n(x) = \tr(\phi(x)S),\] 
and the continuity (actually the analyticity) of $\cL(\nu)$ on $]a,b[$ follows from 
\cite[Cor.~V.4.4]{Ne00}. 
\end{prf}

The preceding theorem generalizes in an obvious way to $\Sesq(\cV)$-valued functions, where 
the corresponding measure $\mu$ has avalues in the cone $\Sesq(\cV)_+$. One can use the 
same arguments as in the proof of Bochner's Theorem.

The following lemma sharpens the ``technical lemma'' in \cite[App.~A]{KL83}. 

\begin{lem} \mlabel{lem:a.2} 
Let $U_t = e^{itH}$ be a unitary one-parameter group on 
$\cH$, $E$ the corresponding spectral measure, $v \in \cH$, 
$E^v := \la E(\cdot)v,v\ra$, $\beta > 0$ and 
$\phi(t) := \la U_t v, v\ra  = \int_\R e^{it\lambda}\, dE^v(\lambda).$ 
Then the following are equivalent: 
\begin{itemize}
\item[\rm(i)] There exists a continuous function $\psi$ on $\oline{\cD_\beta}$, 
holomorphic on $\cD_\beta$, such that $\psi\res_{\R} = \phi$.
\item[\rm(ii)] $\cL(E^v)(\beta) = \int_\R e^{-\beta \lambda}\, dE^v(\lambda) < \infty$. 
\item[\rm(iii)] $v \in \cD(e^{-\frac{\beta}{2}H}).$
\end{itemize}
\end{lem}

\begin{prf} That (i) implies (ii) follows from \cite[p.~311]{Ri66}. 
If, conversely, (ii) is satisfied, then 
$\psi(z) := \cL(E^v)(-iz)$ is defined on $\oline{\cD_\beta}$, holomorphic on $\cD_\beta$ 
and $\psi\res_\R = \phi$. 
Finally, the equivalence of (ii) and (iii) follows immediately from the 
definition of $e^{-\frac{\beta}{2}H}$ in terms of the spectral measure~$E$. 
\end{prf}

\begin{lem} \mlabel{lem:b.4} {\rm(Criterion for the existence of $\cL(\mu)(x)$)} 
Let $\cV$ be a Hilbert space and $\mu$ be a finite 
$\Herm(\cV)_+$-valued Borel measure on $\R$, so that we can consider 
its Laplace transform $\cL(\mu)$, taking values in $\Herm(\cV)$, whenever 
the integral 
\[ \tr(\cL(\mu)(x)S) = \int_\R e^{-\lambda x}\, d\mu^S(\lambda) \quad \mbox { for } \quad 
d\mu^S(\lambda) = \tr(d\mu(\lambda)S),\] 
exists for every positive trace class operator $S$ on $\cV$.
This is equivalent to the finiteness of the integrals 
$\cL(\mu^v)(x)$ for every $v \in \cV$, where $d\mu^v(\lambda) = \la d\mu(\lambda)v,v\ra$. 
\end{lem}

\begin{prf} For $x \in E$, the existence of $\cL(\mu)(x)$ implies the finiteness of the 
integrals $\cL(\mu^v)(x)$ for $v \in \cV$. Suppose, conversely, that all these 
integrals are finite. Then we obtain by polarization a hermitian form 
\[ \beta(v,w) := \int_\R e^{-\lambda x}\, \la d\mu(\lambda)v,w \ra \]  
on $\cV$. We claim that $\beta$ is continuous. 
As $\cV$ is in particular a Fr\'echet space, it suffices to show that, for every 
$w \in \cV$, the linear functional $\lambda(v) := \beta(v,w)$ is continuous 
(\cite[Thm.~2.17]{Ru73}). 

The linear functionals $\lambda_n(v) := \int_{-n}^n e^{-\lambda x}\, \la d\mu(\lambda)v,w \ra$ 
are continuous because $\mu$ is a bounded measure and the functions 
$e_x$ are bounded on bounded intervals. By the Monotone Convergence 
Theorem, combined with the Polarization Identity,
 $\lambda_n \to \lambda$ holds pointwise on $\cV$, and this implies the continuity 
of $\lambda$ (\cite[Thm.~2.8]{Ru73}). 

For a positive trace class operators $S = \sum_n \la \cdot, v_n \ra v_n$ with 
$\tr S = \sum_n \|v_n\|^2 < \infty$, we now obtain 
\[ \cL(\mu^S)(x) 
= \sum_n \cL(\mu^{v_n})(x)
= \sum_n \beta(v_n,v_n) 
\leq \|\beta\|  \sum_n \|v_n\|^2 < \infty.\qedhere\] 
\end{prf}

\end{document}